\documentclass{amsart}
\usepackage{amsmath, amssymb, euscript,  enumerate}

\newtheorem{thm}{Theorem}[section]

\newtheorem{lem}[thm]{Lemma}
\newtheorem{prop}[thm]{Proposition}
\newtheorem{claim}[thm]{Claim}

\theoremstyle{remark}

\numberwithin{equation}{section}

\newcommand{\R}{\mathbb R}

\newcommand{\C}{\mathcal C}

\newcommand{\D}{\nabla}
\newcommand{\intr}{\int_{B_\rho}}
\newcommand{\eu}{\mathrm{eucl}}
\newcommand{\diam}{\mathrm{diam}}
\newcommand{\Int}{\mathrm{Int}}
\newcommand{\ir}{\mathrm{injrad}}
\newcommand{\vol}{\mathrm{Vol}}

\newcommand{\dist}{\mathrm{dist}}

\def\XXint#1#2#3{{\setbox0=\hbox{$#1{#2#3}{\intr}$}
     \vcenter{\hbox{$#2#3$}}\kern-.5\wd0}}

\begin{document}

\title{Eternal Solutions to the Ricci Flow on $\R^2$}

 \author{P. Daskalopoulos$^*$}
\address{Department of Mathematics, Columbia University, New York, USA}
\email{pdaskalo@math.columbia.edu}
\thanks{$*:$ Partially supported by the NSF grants DMS-01-02252, DMS-03-54639  and the EPSRC in the UK}
 \author{N. Sesum}
\address{Department of Mathematics, Columbia University, New York, USA}
\email{natasas@math.columbia.edu}

\begin{abstract} 
We provide the classification  of eternal (or  ancient) solutions of the two-dimensional 
Ricci flow, which is equivalent to the fast diffusion equation $ \frac{\partial u}{\partial t}  = \Delta \log u $
on $ \R^2  \times \R.$ We show that,  under the necessary assumption that for every $t \in \R$, the solution  $u(\cdot, t)$ defines a complete metric of bounded curvature  and  bounded width, $u$   is a gradient soliton of the form $
  U(x,t) = \frac{2}{\beta \, ( |x-x_0|^2 + \delta  \, e^{2\beta t})}$, 
for some $x_0 \in \R^2$ and some  constants $\beta >0$ and $\delta >0$.

\end{abstract}

\maketitle
 
\section{Introduction}\label{sec-0}

We consider  {\em eternal} solutions of the logarithmic fast diffusion  equation 
\begin{equation} \label{eqn-u}
\frac{\partial u}{\partial t}  = \Delta \log u  \qquad \mbox{on} \,\,    \R^2  \times \R.
\end{equation}
This equation represents   the evolution of the conformally flat   
metric  $g_{ij} = u\, I_{ij}$  under the {\it Ricci Flow}
$$
\frac{\partial g_{ij}}{\partial t} = -2 \, R_{ij}.
$$
The equivalence follows from the observation that the metric 
 $g_{ij} = u\, I_{ij}$ has scalar curvature $R = - (\Delta \log u) /u$
 and in two dimensions $R_{ij} = \frac 12 \, R\, g_{ij}$.
 

Equation \eqref{eqn-u} arises also in physical applications, as a
model for long Van-der-Wals interactions in thin  films of a fluid
spreading on a solid surface, if certain nonlinear fourth order
effects are neglected, see \cite{dG,B,BP}.

Our goal in this paper is to provide the  classification of  eternal solutions of equation \eqref{eqn-u} 
under the assumption that for every $t \in \R$, the solution  $u(\cdot, t)$ defines a complete metric of bounded curvature  and  bounded width.  This result is essential in establishing the type II  collapsing  of complete (maximal) solutions of the Ricci flow  \eqref{eqn-u} on  $\R^2 \times [0,T)$ 
with finite area $\int_{\R^2} u(x,t) \, dx < \infty$  (c.f.  in \cite{DD2} for the result in the radially symmetric case).  For an extensive list of  results on  the   Cauchy  problem $u_t=\Delta \log u$ on  $\R^2 \times [0,T)$, $u(x,0)=f(x)$ we refer the reader  to \cite{DD1}, \cite{ERV}, \cite{Hsu}  and  \cite{PT}.

In \cite{DH}   we introduced   the width $w$ of the metric $g=u\, I_{ij}$.
Let $F:\R^2\to[0,\infty)$ denote a 
proper function $F$, such that    $F^{-1}(a)$ is compact for every
$a\in [0,\infty)$.  The width of $F$ is defined to be the supremum
of the lengths of the level curves of $F$, namely 
$w(F) = \sup_c L\{F=c\}.$
The width $w$ of the metric $g$ is defined  to be the infimum
$$w(g) = \inf_F w(F).$$

We will assume, throughout this paper, that $u$ is smooth, strictly positive and satisfies the following conditions:

\noindent The width of the metric $g(t)=u(\cdot,t)\, I_{ij}$ is finite, namely
\begin{equation}\label{eqn-pw}
w(g(t)) < \infty, \qquad \forall t\in \R.
\end{equation}

\noindent The scalar curvature
$R$ satisfies the $L^\infty$-bound
\begin{equation}\label{eqn-bc}
\| R(\cdot,t) \|_{L^\infty(\R^2)} < \infty, \qquad \forall  t \in \R.
\end{equation}

Our goal is to prove the following classification result.

\begin{thm}\label{thm-eternal}
Assume that $u$ is a positive smooth   eternal solution of equation
\eqref{eqn-u}  which defines a complete metric and satisfies conditions
\eqref{eqn-pw}-\eqref{eqn-bc}. Then, $u$ is a gradient soliton of
the form 
\begin{equation}\label{eqn-soliton}
  U(x,t) = \frac{2}{\beta \, ( |x-x_0|^2 + \delta  \, e^{2\beta t})}
\end{equation}
for some $x_0 \in \R^2$ and some  constants $\beta >0$ and $\delta >0$. 
\end{thm}

Under the additional assumptions that the scalar
curvature $R$ is globally bounded on $\R^2 \times \R$ and assumes its maximum at an interior point $(x_0,t_0)$, with
 $-\infty < t_0 < +\infty$, i.e., $R(x_0,t_0) = \max_{(x,t) \in  \R^2 \times \R} R(x,t)$, 
Theorem \ref{thm-eternal} follows from  the result of   R. Hamilton 
on eternal solutions of the Ricci Flow in \cite{H}. 
However, since in general $\partial R/\partial t \geq 0$, without
 this rather  restrictive assumption on the maximum curvature, Hamilton's result does not apply. 

Before we begin with the proof of Theorem \ref{thm-eternal}, let us
give a few remarks.

\smallskip

\noindent{\bf Remarks:} 
\begin{enumerate}[(i)]

\item The bounded width assumption \eqref{eqn-pw} is necessary. If this condition is not satisfied, then \eqref{eqn-u} admits  other solutions,  in particular  the flat (constant) solutions. 

\item It is shown in  \cite{DH} that maximal solutions $u$ of the initial value problem
$u_t=\Delta \log u$ on  $\R^2 \times [0,T)$,  $u(x,0)=f(x)$ which vanish at time $T < \infty$   satisfy the width bound $ c\, (T-t) \leq w(g(t))  \leq C\, (T-t)$ and the maximum curvature bound 
$c\, (T-t)^{-2}  \leq R_{\max}(t)   \leq C\, (T-t)^{-2} $ for some constants $c >0$ and $C< \infty$, independent of $t$. Hence, one may  rescale $u$ near $t \to T$ and pass to the limit to
obtain an eternal solution of equation $\eqref{eqn-u}$ which satisfies the bounds \eqref{eqn-pw} and
\eqref{eqn-bc} (c.f.  in \cite{DD2} for the  radially symmetric case). Theorem \ref{thm-eternal} provides then a classification of the limiting solutions.

\item Since $u$ is strictly positive at all $t < \infty$,  it follows
that $u(\cdot,t)$ must have infinite area, i.e.,
\begin{equation}\label{eqn-ina}
\int_{\R^2} u(x,t) \, dx = +\infty, \qquad \forall t \in \R.
\end{equation}
Otherwise, if $\int_{\R^2} u(x,t) \, dx < \infty$,  for some $t<
\infty$, then by the results in \cite{DD1} the  solution $u$ must
vanish at time $t+T$, with $T= (1/4\pi) \int_{\R^2}  u(x,t) \, dx$,
or before.
\item The proof of Theorem \ref{thm-eternal} only uses that actually
$u$ is an {\em ancient} solution of equation $u_t = \Delta \log u$
on $\R^2 \times (-\infty, T)$, for some $T < \infty$, such that
$\int_{\R^2} u(x,t) \, dx = + \infty$. 
\item Any eternal solution of equation \eqref{eqn-u} satisfies
$u_t \geq 0$.  This is an immediate consequence of the Aronson-B\'enilan inequality (or the maximum
principle on $R=-u_t/u$), which in the case of a solution  on $\R^2 \times [\tau,t)$
states as $u_t \leq u/(t-\tau)$. Letting, $\tau \to -\infty$,  we obtain 
for an eternal  solution  the  bound  $u_t \leq 0$. 
\item Any eternal solution of equation \eqref{eqn-u} satisfies
$R >0$.  By the previous remark, $R= - u_t/u \geq 0$. Since, $R$ evolves
by $R_t = \Delta_g R + R^2$ the strong maximum principle guarantees that $R >0$ or $R \equiv 0$ at all times. Solutions 
with $R \equiv 0$ (flat) violate condition \eqref{eqn-pw}. Hence, $R >0$ on $ \R^2 \times \R$.  
\end{enumerate}

\medskip

\noindent{\bf Acknowledgments.} We are grateful to  S. Brendle, B. Chow, R.
Hamilton, L. Ni, Rafe Mazzeo and F. Pacard for  enlightening discussions in the course of this
work. This paper  was completed while the first author was a Visiting   Fellow  in  the Department Mathematics, Imperial College. She is grateful to this institution for its hospitality and support.

\section{A priori estimates}\label{sec-1}

 We will establish  in this section the asymptotic behavior,  as $|x| \to \infty$, for any  eternal solution  of equation \eqref{eqn-u} which  satisfies  the conditions  \eqref{eqn-pw}-\eqref{eqn-bc}. We will  show that
there exists constants $c(t) >0$ and $C(t) < \infty$ such that 
\begin{equation}\label{asb-u} 
c(t) \, |x|^{-2} \leq  u(x,t) \leq C(t) \,  |x|^{-2}, \qquad t \in \R.
\end{equation}
This bound  is crucial  in the proof of Theorem \ref{thm-eternal}.

We begin with the following lower bound which is a consequence of the results  in \cite{PT}.

\begin{prop}\label{prop-max}
Assume that $u$ is a positive smooth   eternal solution of equation
\eqref{eqn-u}  which defines a complete metric and satisfies condition \eqref{eqn-bc}. Then, 
\begin{equation}\label{equation-decay}
\frac 1{u(x,t)} \leq  O \left (r^2 \, \log^2 r \right ), \qquad \mbox{as} \,\,\, r=|x|\to\infty, \quad \forall t  \in \R.
\end{equation}
\end{prop} 
\begin{proof}
For $t_0 \in \R$  fixed, let $\bar u$ denote the maximal solution of the Cauchy problem 
\begin{equation}
\begin{cases}
\label{equation-comp}
 \bar u_t  = \Delta\log \bar u,  \qquad  &\mbox{on} \,\, \R^2 \times (0, \infty) \cr \bar u(x,0) =  u(x,t_0) \qquad & x\in \R^2.
\end{cases}
\end{equation}
It follows by the results of Rodriguez, Vazquez and Esteban in \cite{ERV} that $u$ satisfies the growth condition
\begin{equation}\label{eqn-max}
\frac{1}{\bar u(x,t)} \le O \left (\frac{r^2 \, \log^2 r }{2t} \right ), \qquad \mbox{as} \,\,\, r=|x|\to\infty, \quad \forall t >0.
\end{equation}
In particular, $\bar u$ defines a complete metric. Let us denote by $\bar R$ the
curvature of the metric $\bar g(t) = \bar u \, I_{ij}$.

\begin{claim}
\label{claim2.2}  
There exists $\tau >0$ for which
$\sup_{\R^2}|\bar R(\cdot,t))| \le C(\tau)$ for $t \in [0,\tau].$
\end{claim}

>From the claim the proof of the Proposition readily
follows by the uniqueness result of Chen and Zhu (\cite{ChZh}; see
also \cite{PT}). Indeed, since both $u$ and $\bar u$ define complete
metrics with bounded curvature, the uniqueness result in \cite{ChZh}
implies that $\bar u(x,t)=u(x,t+t_0)$, for $t \in [0,\tau]$. Hence,
$u$ satisfies \eqref{eqn-max}
 which readily implies
\eqref{equation-decay}, since $u$ is decreasing in
 time. 
 
\noindent{\em Proof of Claim \ref{claim2.2}.}
Since $\bar g(0) = g(t_0)$ and 
$$0 < \bar R(\cdot,  0) = R(\cdot,t_0) \leq C_0$$ the  classical result of Klingenberg (see \cite{GHL}) implies  the injectivity radius 
bound  
$r_0 = \ir(\R^2, \bar g(0)) \ge \frac{\pi}{C_0}$. Moreover,
\begin{equation}
\label{equation-volume}
\vol_{\bar g(0)}B_{\bar g(0)}(x,r) \ge V_{C_0}(r),
\end{equation}
where $V_{C_0}(r)$ is the   volume 
of a ball of radius $r$ in a space form of constant sectional curvature 
$C_0$ (see \cite{J} for (\ref{equation-volume})). 
We will prove the desired  curvature bound  using
(\ref{equation-volume}) and Perelman's pseudolocality theorem
(Theorem $10.3$ in \cite{Pe}). Let $\epsilon, \delta > 0$ be as in the 
pseudolocality theorem.  Choose $r_1 < r_0$ such that for all $r\le r_1$,
we have $V_{C_0}(r) \ge (1-\delta) \, r^2$. By (\ref{equation-volume}), for all
$x_0\in \R^2$ we have $\vol_{\bar g (0)}B_{\bar{g}(0)}(x_0,r) \ge (1-\delta)r^2$.
Since also $\sup_{\R^2}|\bar R(\cdot,0)| \le C_0$,  the pseudolocality theorem implies
the bound
$$|\bar R(x,t)| \le (\epsilon r_1)^{-2} \qquad \mbox{whenever}
\,\,\,0 \le t \le  (\epsilon r_1)^2, \,\,\, \dist_t(x,x_0) < \epsilon r_1.$$
Since the previous estimate does not depend on $x$, we obtain the uniform bound 
$$\sup_{\R^2}|\bar R(\cdot,t)|  \le (\epsilon r_1)^{-2} \qquad  \forall t\in [0,(\epsilon r_1)^2]$$
finishing the proof of the claim. 
\end{proof}

We will next perform the  cylindrical change of coordinates, setting
\begin{equation}\label{eqn-dv}
v(s,\theta,t) = r^2 \, u(r,\theta,t), \qquad s=\log r
\end{equation}
where $(r,\theta)$ denote polar coordinates. It is then easy to see
that the function $v$ satisfies the equation
\begin{equation}\label{eqn-v}
v_t =  (\log v)_{ss} + (\log v)_{\theta\theta}, 
 \qquad \mbox{for} \,\, 
(s,\theta,t) \in {\mathcal C}_{\infty} \times \R
\end{equation}
with $\C_{\infty}$ denoting the infinite cylinder $\C_\infty= \R \times [0,2 \pi]$. 
Notice that the nonnegative curvature condition $R \geq 0$ implies that 
\begin{equation}\label{eqn-v2}
\Delta_c \log v := (\log v)_{ss} + (\log v)_{\theta\theta} \leq 0
\end{equation}
namely that $\log v$ is superharmonic in the cylindrical $(s,\theta)$ 
coordinates.
Estimate \eqref{asb-u} is equivalent to:

\begin{lem}\label{lem-asv} For every $t \in R$, there exist constants $c(t) >0$ and $C(t) <\infty$ such that
\begin{equation}\label{eqn-asv}
c(t) \leq  v(s,\theta,t) \leq  C(t), \qquad  (s,\theta)  \in \C_\infty.  
\end{equation}
\end{lem}

The proof of Lemma \ref{lem-asv} will be done in several steps. We will first establish the bound from below which only uses that the curvature  $R \geq 0$ and that the metric is complete.

\begin{prop}\label{claim-lower-bar-v} If $u(x,t)$ is a maximal solution of (\ref{eqn-u}) that defines
a metric of positive curvature, then for every $t \in R$, there exists a  constant $c(t) >0$ such that
\begin{equation}\label{eqn-asv2}
v(s,\theta,t) \geq  c(t), \qquad  (s,\theta)  \in \C_\infty.  
\end{equation}
\end{prop}
\begin{proof}
Fix a $t \in \R$. We will show that $(\log v)^-= \max (-\log v,0)$ is bounded above.
We begin by observing that 
$$\Delta_c (\log v)^- =  (\log v)^-_{ss} +  (\log v)^-_{\theta\theta} \geq 0$$
since $\Delta_c \log v =  -R\, v \leq 0$. Hence, setting 
$$V^-(s,t) = \int_0^{2\pi} (\log v)^-(s,\theta,t) \, d\theta$$
the function $V^-$ satisfies $V^-_{ss} \geq 0$, i.e., $V_s^-$ is increasing
in $s$. It follows that the limit $\gamma = \lim_{s\to\infty}V^-_s(s,t)  \in [0, \infty]$ exists. If $\gamma >0$,  then 
$V^-(s,t)  >  \gamma_1 \, s$,
with $\gamma_1=\gamma/2$  for $s \geq s_0$ sufficiently large,
which contradicts the pointwise bound \eqref{equation-decay}, which
when expressed in  terms of $v$ gives $v(s,t) \leq c(t)/s^2$, for $s \geq s_0$ sufficiently large. We conclude that $\gamma=0$. Since $V^-_s$
is increasing, this implies that $V^-_s \leq 0$, for all $s \geq s_0$, i.e.
$V^-(s,t)$ is decreasing in $s$ and therefore bounded above. 

We will now derive a pointwise bound on $(\log v)^-$. Fix a point
$(\bar{s},\bar{\theta})$ in cylindrical coordinates, with  $\bar{s} \ge
s_0+2\pi$ and let $B_{2\pi} = \{(s,\theta):  |s-\bar{s}|^2 +
|\theta-\bar{\theta}|^2 \le (2\pi)^2\}$. 
By the mean value inequality for sub-harmonic functions, we obtain 
$$(\log v)^-(\bar{s},\bar{\theta},t) \leq \frac{1}{|B_{2\pi}|}\int_{B_{2\pi}}(\log v)^- (s,\theta,t) \, ds \, d\theta.$$
Since $B_{2\pi} \subset Q = \{|s-\bar{s}| \le 2\pi, \,\, \theta \in
[-4\pi,4\pi]\}$, we conclude the bound
$$ (\log v)^-(\bar{s},\bar{\theta},t) \leq C \int_Q (\log v)^- (s,\theta,t) \, ds \, d\theta \leq 4 C \int_{\bar s-2\pi}^{\bar s+2\pi} V^-(s,t)\, ds \leq \tilde C$$ finishing the proof.
\end{proof}

The estimate from above on $v$ will be based on the following 
integral bound. 

\begin{prop} \label{prop-log+}
Under the assumptions of Theorem \ref{thm-eternal},
for every $t \in \R$, we have
\begin{equation}
\label{eqn-blog+}
\sup_{s\ge 0}\int_0^{2\pi}(\log v(s,\theta,t))^+ d\theta < \infty. 
\end{equation}
\end{prop}

\begin{proof}
Fix $t \in \R$ 
and define 
$$V(s,t) = \int_0^{2\pi}\log v(s,\theta,t)  \, d\theta.$$ 
Since $\Delta_c \log v \leq 0$, $V$ satisfies $V_{ss} \le 0$, i.e. 
$V_s$ decreases in $s$.  Set 
$\gamma = \lim_{s\to\infty}V_s \in [-\infty,\infty)$.

\begin{claim}
\label{lem-zero}
We have  $\gamma = 0$, i.e.  $V$ is increasing in $s$, and
$\lim_{s \to \infty} V(s,t) < \infty$. 
\end{claim}

\noindent{\em Proof of Claim:} 
We  consider the following two  cases:

\noindent{\em Case 1:  $\gamma \ge 0$.} Then since $V_s$ decreases in $s$, we have $V_s \ge
0$ for all $s$, which implies that  there is $\beta =
\lim_{s\to\infty}V(s,t) \in (-\infty,\infty]$.
If $\gamma > 0$, there is $s_0$ such that for $s\ge s_0$,
\begin{equation}
\label{equation-lim-v}
\int_0^{2\pi}\log v(s,\theta,t) \, d\theta \ge \gamma_1 s,
\end{equation}
where $\gamma_1 = \gamma/2$ . 
\begin{enumerate}
\item[(1a)]
If $\beta < \infty$, then for $s >> 1$,
$$V(s,t_0) \le \beta,$$
which contradicts (\ref{equation-lim-v}) for big $s$, unless $\gamma = 0$.
\item[(1b)] 
If $\beta = \infty$ and $\gamma > 0$ we will   derive a
contradiction using the boundness of the width.  The function $\log
v$ satisfies 
$\Delta_c \log v = -R \, v \le 0$
that is, $\log v$ is the superharmonic function.  Fix a point
$(\bar{s},\bar{\theta})$ in cylindrical coordinates, with  $\bar{s} \ge
s_0$ and let $B_{2\pi} = \{(s,\theta):  |s-\bar{s}|^2 +
|\theta-\bar{\theta}|^2 \le (2\pi)^2\}$. 
By the mean value inequality for superharmonic functions, we obtain 
$$\log v(\bar{s},\bar{\theta},t) \ge  \frac{1}{|B_{2\pi}|}\int_{B_{2\pi}}\log v (s,\theta,t) \, ds \, d\theta.$$
Expressing  $\log v = (\log v)^+ - (\log v)^-$, we observe that 
$$(\log v)^-(s,\theta,t) \le C(t) + 2\log s, \qquad \mbox{on}  \,\, B_{2\pi}$$
from the bound \eqref{equation-decay}. 
Since 
$Q = \{|s-\bar{s}|\le 1, \,\, |\theta-\bar{\theta}|\le \pi\}$ is
contained in $B_{2\pi}$, we obtain 
\begin{equation*}
\begin{split}
\int_{B_{2\pi}} (\log v)^+  &\geq  \int_{\bar{s}-1}^{\bar{s}+1}\int_{\bar{\theta}-\pi}^{\bar{\theta}+\pi}(\log v)^+ 
=\int_{\bar{s}-1}^{\bar{s}+1}\int_0^{2\pi} (\log v)^+ \\
&\geq V(\bar s,t) - C(t) + 2\log s.
\end{split}
\end{equation*}
Combining the above  with \eqref{equation-lim-v} yields the bound 
\begin{equation}\label{eqn-plogv}
\log v(\bar{s},\bar{\theta},t) \ge   \gamma_1 s - 2C(t) - 4\log s \ge \bar{\gamma}s
\end{equation}
for $\bar{s} >> 1$ and $\bar{\gamma} < \gamma_1$.
We conclude that 
\begin{equation}
\label{equation-contra-width}
v(\bar{s},\bar{\theta},t) \ge e^{\bar{\gamma}\bar{s}} \to \infty, \qquad \mbox{as} \,\,
\bar{s}\to\infty.
\end{equation}
We will next show that (\ref{equation-contra-width}) actually  contradicts  our
bounded    width  condition $\bar w(\bar g(t_0)) < \infty$. Indeed, let $F:\mathrm{R}^2 \to
[0,\infty)$ be a proper function on the plane. Denoting by
$L_g\{F=c\}$ the length of the $c$-level curve $\sigma_c$ of $F$,  measured with respect to metric $g(t)=u(\cdot,t)\, I_{ij}$, and using the bound \eqref{eqn-plogv}, namely  $ u(x,t)|x|^2 = v(\log |x|,\theta,t_0) \ge e^{\bar{\gamma}\log|x|}$, we obtain 
\begin{eqnarray*}
L_g(\sigma_c) &=& \int_{\sigma_c}\sqrt{u} \, d\sigma_c \\
&\ge& e^{\min_{x\in \sigma_c}((\bar{\gamma}/2-1)\log|x|)}L_{\eu}(\sigma_c).
\end{eqnarray*}
If $\bar{\gamma} \ge 2$, then  denoting  by $|x_c| = \min_{x\in\sigma_c}|x|$ we have 
\begin{eqnarray*}
L_g(\sigma_c) &\ge& e^{(\bar{\gamma}/2-1)\log|x_c|}L_{\eu}(\sigma_c) \ge e^{(\bar{\gamma}/2-1)\log|x_c|}2\pi|x_c| \\
&=& |x_c|^{\bar{\gamma}/2} \to \infty, \qquad \mbox{as} \,\,\, c\to\infty
\end{eqnarray*}
where we have used the fact that the euclidean circle centred at the origin  of
radius $|x_c|$ is contained in the region bounded by the curve $\sigma_c$.\\
If $0 < \bar{\gamma} < 2$, then $\alpha = 1 - \bar{\gamma}/2 > 0$ and
since we may assume that the origin is contained in  the interior of the region bounded by the
level curve $\sigma_c$, denoting  by $|x_c| = \max_{x\in \sigma_c}|x|$, 
we obtain 
\begin{eqnarray*}
L_g(\sigma_c) &\ge& e^{-\alpha\log|x_c|}L_{\eu}(\sigma_c) \\
&\ge& \frac{1}{|x_c|^{\alpha}}|x_c| \\
&=& |x_c|^{1-\alpha} \to \infty \qquad \mbox{as} \,\,\, c\to\infty.
\end{eqnarray*}
Here we have used the fact that 
$$L_{\eu}(\sigma_c) \ge \diam_{\eu}(\overline{\Int{\sigma_c}}) \ge \max_{x\in\sigma_c}|x| = |x_c|.$$
We conclude that  
$$\sup_c L_g\{F=c\} \ge M$$
for any proper function $F$, which implies that 
$w(g(t))  \ge M$, 
contradicting our width bound \eqref{eqn-pw}. 
\item[(1c)]
If $\beta = \infty$ and $\gamma = 0$, we can argue similarly as in (1b),
with the only difference that now  $V(s,t)$ increases in $s$
and $V(s,t) >> M$ for $s >> 1$, where $M$ is an arbitrarily large constant.
The mean value property applied to $\log v$ now shows   that
$\log v$ can be made arbitrarily large for $s >> 1$. The rest of the argument is the same as in (1b).
\end{enumerate}
So far we have proven  that if $\gamma \ge 0$, the only posibility is that $\gamma = 0$ and $\beta < \infty$.

\noindent{\em Case 2:}  $\gamma < 0$. Then,  there is $s_0$ such that for $s\ge s_0$,
\begin{equation}
\label{equation-small-v}
\int_0^{2\pi}\log v(s,\theta,t) d\theta < -|\gamma|s.
\end{equation}
On the other hand, by (\ref{equation-decay}) we have
$$V(s,\theta,t) \ge \frac{C(t)}{s^2}, \qquad \mbox{for} \,\, s >> 1$$ which implies
\begin{eqnarray*}
-|\gamma|s &>& \int_0^{2\pi}\log v \, d\theta \\
&\ge& 2\pi\log(C(t)) - 4\pi\log s = \tilde{C}(t) - 4\pi\log s,
\end{eqnarray*}
which is not possible for large  $s$. 

We have shown that the only possibility is that $\gamma = 0$
(which means $V(s,t)$ is increasing in $s$) and $\beta < \infty$. 
This finishes the proof of the Claim. 

We can now finish the proof of the Proposition. Expressing 
$$\int_0^{2\pi}(\log v)^+(s,\theta,t) \, d\theta   = V(s,t)  + 
\int_0^{2\pi}(\log v)^-(s,\theta,t) \, d\theta$$
and using  that $V (s,t) \le \beta$ (as shown in Claim \ref{lem-zero}) together with Proposition \ref{claim-lower-bar-v}, we readily obtain 
$$\int_0^{2\pi}(\log v)^+ (s,\theta,t) \, d\theta\le \beta + 2\pi\delta(t) < \infty$$
as desired. 
\end{proof}

To prove Lemma \ref{lem-asv}, it only remains to show that
$\|v(\cdot,t) \|_{L^\infty} < \infty$. The proof of this bound will
use the ideas  of Brezis and Merle (\cite{br-mer}), 
including the following result which we state for the reader's convenience. 

\begin{thm}[Brezis-Merle]
\label{theorem-br-mer}
Assume $\Omega \subset  \mathrm{R}^2$ is a bounded domain and let $u$ be a 
solution of 
\begin{equation}\label{eqn-ell}
\begin{cases}
-\Delta u = f(x) \qquad &\mbox{in} \,\,\, \Omega, \cr
\,\,\,\, \quad u = 0 \qquad &\mbox{on} \,\,\, \partial\Omega,
\end{cases}
\end{equation}
with $f\in L^1(\Omega)$. Then for every $\delta \in (0,4\pi)$ we have
$$\int_{\Omega}e^{\frac{(4\pi-\delta)|u(x)|}{||f||_{L^1(\Omega)}}}dx
\le \frac{4\pi^2}{\delta}(\diam\Omega)^2.$$
\end{thm}

\medskip

\noindent{\em Proof of Lemma \ref{lem-asv}.}  The proof follows the ideas of Brezis and Merle in \cite{br-mer}. Set  $w=\log v $ so that 
$$\Delta w = -R \, e^w$$
with $R$ denoting the scalar curvature. 

Fix $\epsilon \in (0,1)$. Since $\int_{\mathrm{R}^2}R \, e^w < \infty$, there is
$s_0$ such that for all $s\ge s_0$
$$\int_{s-2\pi}^{s+2\pi}\int_0^{2\pi}R \, e^w\, d\theta \, d\rho < \epsilon.$$
Fix $\bar s > s_0$ and set  $B_{2\pi}(\bar{s},\bar{\theta}) = \{(s,\theta)|\,\, (s-\bar{s})^2 +
(\theta-\bar{\theta})^2 \le (2\pi)^2\} \subset \{(s,\theta)|\,\,
|s-\bar{s}| \le 2\pi, \,\,\, |\theta-\bar{\theta}| \le 2\pi\} =
Q(\bar{s},\bar{\theta})$ so that 
$$\int_{B_{2\pi}(\bar{s},\bar{\theta})}  R \, e^w < \epsilon.$$
We shall denote by $B_{2\pi}$  the  ball
$B_{2\pi}(\bar{s},\bar{\theta})$ and by $Q$  the  cube $Q(\bar{s},\bar{\theta})$.
Let $w_1$ solve  problem \eqref{eqn-ell} on $\Omega=B_{2\pi}$ with $f=Re^w \in L^1(B_{2\pi})$ and $\|f\|_{L^1(B_{2\pi})} < \epsilon$. 
By Theorem \ref{theorem-br-mer} there exists a constant $C >0$ 
for which 
$$\int_{B_{2\pi}}e^{\frac{(4\pi-\delta)|w_1(x)|}{||f||_1}} \le
\frac{C}{\delta}.$$
Taking  $\delta = 4\pi - 1$ we obtain 
\begin{equation}
\label{equation-w1-p}
\int_{B_{2\pi}}e^{\frac{(4\pi-\delta)|w_1(x)|}{\epsilon}} \le C.
\end{equation}
Combining (\ref{equation-w1-p}) and Jensen's inequality gives the estimate 
\begin{equation}
\label{equation-w1}
||w_1||_{L^1(B_{2\pi})} \le \tilde{C}(\epsilon).
\end{equation}
The difference  $w_2 = w-w_1$ satisfies  $\Delta w_2 = 0$ on $B_{2\pi}$. Hence by the  mean value
inequality 
\begin{equation}
\label{equation-w2}
||w_2^+||_{L^{\infty}(B_{\pi})} \le C||w_2^+||_{L^1(B_{2\pi})}.
\end{equation}
Since
$w_2^+ \le w^+ + |w_1|$ 
combining (\ref{equation-w1}) and (\ref{eqn-blog+}) yields
$$||w_2^+||_{L^1(B_{2\pi})} \le C.$$
Expressing  $Re^w = R\, e^{w_2}e^{w_1}$ and observing that 
\begin{itemize}
\item
$R$ is bounded and  (\ref{equation-w2}) holds on $B_{\pi}$
\item
$e^{w_1} \in L^{\frac{1}{\epsilon}}(B_{2\pi})$, where $\epsilon$ can be chosen small, so that $\frac{1}{\epsilon} > 1$, and
\item
$-\Delta w = R e^w \le C \, e^{w_2}$
\end{itemize}
by standard elliptic estimates we obtain the  bound 
\begin{equation}
\label{equation-upper}
||w^+||_{L^{\infty}(B_{\pi/2})} \le C\, ||w^+||_{L^1(B_{\pi})} + C\, ||e^{w_1}||_{L^p(B_{\pi})} 
\le \tilde{C},
\end{equation}
where $p > 1$, finishing the proof of the Lemma.
\qed

Fix $\tau \in \R$ and define the inifinite cylinder
$$Q_+(\tau) = \{ (s,\theta,t): \, s \geq 0, \theta \in [0,2\pi], \tau -1 \leq t \leq \tau \, \}.$$
By Lemma \ref{eqn-asv} 
\begin{equation}\label{eqn-bv}
 0 < c(\tau) \leq v(s,\theta,t) \leq C(\tau) < \infty, \qquad  \mbox{on} \,\, Q_+(\tau).
 \end{equation}
It follows that equation \eqref{eqn-v} in uniformly parabolic in $Q_+(\tau)$.
Hence, the bounds \eqref{eqn-bv} combined with classical derivative
estimates for uniformly parabolic equations imply the following:

\begin{lem}\label{lem-3} Under the assumptions of Theorem
\ref{thm-eternal} for every $\tau \in \R$, there exists a constant
$C$ such that 
$$|v_s| \leq \frac C{s}, \quad |v_{ss}| \leq \frac C{s^2}, \quad |v_\theta| \leq C, \quad |v_{\theta s}| \leq \frac C{s}, \qquad \mbox{on} \,\, Q_+(\tau).$$
\end{lem}

 \smallskip
 We will next show that the curvature
$R(x,t)= - (\Delta \log u(x,t)) /u(x,t)$ tends 
 to zero, as $|x| \to
\infty$. 
 We begin by reviewing the
 Harnack inequality satisfied
by the curvature $R$, shown by R. Hamilton \cite{H4} and \cite{H1}. 
In the case of eternal  solutions $u$ of \eqref{eqn-u} which define a complete metric,  it
states as
\begin{equation}\label{eqn-harn1}
 \frac{\partial \log R}{\partial t} \geq |D_g \log R|^2.
 \end{equation}
Since,  $D_g R= u^{-1} \, D R$, equivalently, this gives the inequality
\begin{equation}\label{eqn-harn2}
 \frac{\partial  R}{\partial t} \geq \frac{|D R|^2}{R\,u}. 
\end{equation}
Let $(x_1, t_1), (x_2, t_2)$
be any two points in $\R^2 \times \R$, with $t_2 > t_1$. Integrating
\eqref{eqn-harn2} along the path $x(t) = x_1 + \frac{t - t_1}{t_2-t_1} \,  x_2$, also using the bound  $u(x,t) \leq C(t_1) / |x|^2$, holding for all $t \in [t_1,t_2]$ by \eqref{eqn-pw} and the fact that $u_t \leq 0$,  we find  the more standard in PDE Harnack inequality 
\begin{equation}\label{eqn-harn3}
  R(x_2, t_2) \geq R(x_1,t_1) \, e^{-C \, (\frac{|x_2-x_1|^2}{|x_1|^2 \, 
      (t_2 - t_1)})}.
\end{equation}

One may now combine \eqref{eqn-harn3} with Lemmas \ref{lem-asv} and \ref{lem-3}  to conclude the following:

\begin{lem}\label{lem-4}
 Under the assumptions of Theorem \ref{thm-eternal}
we have
$$\lim_{|x| \to \infty} R(x,t) =0, \qquad \forall t \in \R.$$
\end{lem}
\begin{proof}
For any $t \in \R$ and $r >0$, we denote by $\bar R(r,t)$ the spherical average
of $R$, namely  
$$\bar R(r,t) = \frac 1{|\partial B_r|} \int_{\partial B_r} R(x,t) \, d\sigma.$$
We first claim that
\begin{equation}\label{eqn-14}
\lim_{r \to \infty} \bar R(r,t) =0, \qquad \forall t.
\end{equation}
Indeed, using the cylindrical coordinates, introduced previously,
this claim is equivalent to showing that
$$\lim_{s \to \infty} \int_0^{2\pi} - \frac{(\log v)_{ss}(s,\theta,t) + (\log v)_{\theta \theta}(s,\theta,t)}{v} \, d\theta =0$$
which is equivalent, using Lemma \ref{lem-asv}, to showing  that
$$\lim_{s \to \infty} \int_0^{2\pi} - (\log v)_{ss}(s,\theta,t) \, d\theta =0.$$
But this readily follows from Lemmas \ref{lem-asv} and \ref{lem-3}.

Fix $t \in \R$. To prove  that  $\lim_{|x| \to \infty} R(x,t) = 0$, we use the Harnack inequality \eqref{eqn-harn3} to show  that
$$R(x,t) \leq  C\, \inf_{y \in \partial B_r} R(y,t+1) \leq  C\, \bar R(r,t+1), \qquad \forall x \in \partial B_r,  \,\, \forall t \in R$$
and use \eqref{eqn-14}. 
\end{proof}

Combining the above with classical derivative estimates for linear
strictly parabolic equations,  gives the following.

\begin{lem}\label{lem-5}
Under the assumptions of Theorem \ref{thm-eternal} the radial
derivative $R_r$ of the curvature satisfies 
$$\lim_{|x|\to \infty} |x|\, R_r (x, t) =0, \qquad \forall t \in \R.$$
\end{lem}

\begin{proof}
For any $\rho >1$ we set $\tilde R(x,t) = R(\rho \, x,t)$ and we compute
from the evolution equation $R_t = u^{-1} \Delta R + R^2$ of $R$,
that 
$$\tilde R_t = (\rho^2  u)^{-1} \Delta \tilde R + \tilde R^2.$$
For $\tau <T$  consider the cylinder $Q = \{ (r,t): 1/2 \leq |x| \leq 4, \tau - 1 \leq t \leq \tau \, \}$. 
>From  \eqref{asb-u} we have  $ 0 < c(\tau)  \leq \rho^2 \, u(x,t) \leq C(\tau)< \infty$, for all  $x\in Q$, 
hence $\tilde R$ satisfies a uniformly parabolic equation in $Q$. 
Classical derivative estimates then imply that
$$ |(\tilde R)_r(x,t) | \leq C \, \| \tilde R\|_{L^\infty(Q)}$$
for all $1 \leq |x| \leq 2$, $\tau - 1/2 \leq t \leq \tau$, implying in particular that
$$ \rho \, |R_r (x,\tau) | \leq C \, \|  R\|_{L^\infty(Q_\rho)}$$
for all $\rho \leq  |x|  \leq 2\, \rho$, where $Q_\rho =  \{ (x,t): \rho/2 \leq |x| \leq 4\rho, \tau - 1 \leq t \leq \tau\}$. The proof now follows from Lemma \ref{lem-4}. 
\end{proof}

\section{Proof of Theorem \ref{thm-eternal}}

 Most of the computations here are known  in the case that $u\, (dx_1^2 + dx_2^2)$ defines a metric on a compact surface (see for example in \cite{Ch}). However, in  the non-compact case an exact account of the boundary terms at infinity should be made. 

We begin
by integrating the Harnack inequality $R_t \geq  |DR|^2/Ru$   with respect to the measure
$d\mu = u \, dx$. Since the measure $d\mu$ has infinite area,
we will  intergrate over a fixed ball $B_{\rho}$. At the end of the 
proof we will let $\rho \to \infty$.  
 Using also that $R_t = u^{-1} \Delta R + R^2$
we find
$$ \intr \Delta R \, dx  + \intr R^2 \, u\, dx \geq \intr \frac{ |DR|^2}R \, dx $$
and by Green's Theorem  we conclude
\begin{equation}\label{eqn-eter1}
 \intr \frac{ |DR|^2}R \, dx - \intr R^2 \, u dx  \leq \int_{\partial B_\rho} \frac{\partial R}{\partial \nu} d\sigma. 
\end{equation}
  
Next,  following Chow (\cite{Ch}), we consider the vector $X= \nabla R + R\, \nabla f$, where $f=-\log u$ is the potential function (defined up to a constant) of the scalar curvature, since it
 satisfies $\Delta_g f = R$, with $\Delta_g f= u^{-1} \, \Delta f$ 
 denoting the Laplacian with respect to the conformal metric $g= u\, (dx^2 + dy^2)$. As it was observed in \cite{Ch} $X \equiv 0$ on  Ricci solitons, i.e.,  Ricci solitons are gradient solitons in the direction
 of $\nabla_g f$. 
A direct computation shows 
 $$
 \intr \frac{|X|^2}{R}\, dx =   \intr \frac{|D R|^2}{R} \, dx +  2 \,  \intr \D R \cdot \D f \, dx + \intr R \,  |D f|^2 \, dx.
$$
Integration by parts  implies 
$$ \intr \D R \cdot \D f \, dx  = - \intr R\, \Delta f \, dx +  
 \int_{\partial B_\rho} R\,  \frac{\partial f}{\partial n} \, d\sigma = - \intr R^2 \, u\, dx +  
 \int_{\partial B_\rho} R\,  \frac{\partial f}{\partial n} \, d\sigma
$$
since 
$\Delta f = R\, u$.
Hence
\begin{equation}\label{eqn-eter2}
\begin{split}
\intr \frac{|X|^2}{R}\, dx = \intr \frac{|D R|^2}{R} &\, dx - 2\,  \intr R^2 \, u\, dx \\&+ \intr R \,  |D f|^2 \, dx +  2\, \int_{\partial B_\rho} R\,  \frac{\partial f}{\partial n} \, d\sigma. 
\end{split}
\end{equation}
Combining \eqref{eqn-eter1} and \eqref{eqn-eter2} we find that
\begin{equation}\label{eqn-eter3}
\intr \frac{|X|^2}{R}\, dx \leq  - \left ( \intr R^2 \, u\, dx  -  \intr R \,  |D f|^2 \, dx \right) + I_\rho = - M + I_{\rho}
 \end{equation}
where
$$I_{\rho} =  \int_{\partial B_\rho} \frac{\partial R}{\partial n} \, d\sigma + 2\, \int_{\partial B_\rho} R\,  \frac{\partial f}{\partial n} \, d\sigma.$$
Lemmas \ref{lem-3} - \ref{lem-5} readily imply that
\begin{equation}\label{eqn-eter4}
\lim_{\rho \to \infty} I_\rho =0.
\end{equation}

As in  \cite{Ch},  we will show next that $M \geq 0$ and indeed a complete square
which vanishes exactly on Ricci solitons. 
To this end, we  define the matrix
$$M_{ij} = D_{ij} f  + D_i f\, D_j f - \frac 12  ( |D f|^2 + R\, u ) \, I_{ij}$$
with $I_{ij}$ denoting the identity matrix. A direct computation shows
that $M_{ij} = \D_i \D_j f - \frac 12 \Delta_g f \, g_{ij}$, with $\nabla_i$
denoting covariant derivatives. It is well known that the  Ricci
solitons are characterized by the condition $M_{ij} =0$, (see in  \cite{H1}).  

\noindent{\em Claim:}  
\begin{equation}\label{eqn-eter33}
M:=  \intr R^2 \, u\, dx  -  \intr R \,  |D f|^2 \, dx = 2\, \intr |M_{ij}|^2 \frac 1u \, dx + J_\rho \end{equation}
where 
$$\lim_{\rho \to \infty} J_\rho =0.$$
To prove the claim we first observe that since $\Delta f = R u$
$$
 \intr R^2 \, u =   
 \intr \frac{(\Delta f)^2}{u} \,  dx  = \intr D_{ii} f \, D_{jj} f \, \frac 1u\, dx.$$
Integrating by parts and  using again  that $\Delta f = R u$, we find \begin{equation*}
\begin{split}
 \intr D_{ii} f \, D_{jj} &f\, \frac 1u \, dx = - \intr D_{jii} f \, D_j f \frac 1u \, dx  \\&+ \intr \Delta   f \, D_j f \frac {D_j u}{u^2} \, dx \
+ \int_{\partial B_\rho} R  \, \frac{\partial f}{\partial n} \, d\sigma.
\end{split}
\end{equation*}
Integrating by parts once more we find 
\begin{equation*}
\begin{split}
\intr D_{jii} f \, D_j f & \frac 1u \, dx = - \intr |D_{ij} f|^2  \frac 1u \, dx\\
&+ \intr D_{ij} f \, D_j f \frac {D_i u}{u^2} \, dx + \frac 12  \int_{\partial B_\rho}  \frac{\partial (|Df|^2)}{\partial n}  \, \frac 1u \, d\sigma
\end{split}
\end{equation*}
since
$$ \int_{\partial B_\rho} D_{ij} f \, D_j f \, n_i \, \frac 1u \, d\sigma = \frac 12  \int_{\partial B_\rho} 
\frac{\partial (|Df|^2)}{\partial n}  \, \frac 1u \, d\sigma.$$
Combining the above and using that $D f = - u^{-1} Du$ and $\Delta f = R\, u$ we conclude
\begin{equation}\label{eqn-eter44}
 \intr R^2  u\, dx  =  \intr |D_{ij} f|^2  \frac 1u \, dx + \intr D_{ij} f \, D_i f \, D_j f \frac 1u \, dx  - \intr  R\,  |D  f|^2 \,   dx  + J^1_\rho
\end{equation}
where 
$$J^1_\rho =  \int_{\partial B_\rho} R \,  \frac{\partial f}{\partial n} \, d\sigma -  \frac 12  \int_{\partial B_\rho}  \frac{\partial (|Df|^2)}{\partial n}  \, \frac 1u \, d\sigma.$$
Hence
\begin{equation}\label{eqn-eter5}
M   = \intr |D_{ij} f|^2  \frac 1u \, dx + \intr D_{ij} f \, D_i f \, D_j f \frac 1u \, dx  - 2 \, \intr  R\,  |D  f|^2 \,   dx + J^1_\rho.
\end{equation}

We will now intergrate $|M_{ij}|^2$. 
A direct computation and  $\Delta f = R u$ imply 
\begin{equation}\label{eqn-eter6}
\begin{split}
\intr | M_{ij} |^2 \frac 1u \, dx &=   \intr |D_{ij} f|^2  \frac 1u \, dx + 2  \intr D_{ij} f \, D_i f \, D_j f \frac 1u \, dx \\ &-   \intr R \,  |D f|^2 \, dx  +  \frac 12  \intr   |D f|^4 \, \frac 1u  dx  - \frac 12 \intr R^2\, u\, dx.
\end{split}
\end{equation}
Combining \eqref{eqn-eter5} and \eqref{eqn-eter6} we then find
\begin{equation*}
\begin{split}
M - 2\, \intr | M_{ij} |^2 \frac 1u \, dx &= -  \intr |D_{ij} f|^2  \frac 1u \, dx - 3 \intr D_{ij} f \, D_i f \, D_j f \frac 1u \, dx \\&  -   \intr   |D f|^4 \, \frac 1u  dx +   \intr R^2 u\, dx  + J_\rho^1. 
\end{split}
\end{equation*}
Using \eqref{eqn-eter44} we then conclude that
\begin{equation}\label{eqn-eter7}
\begin{split}
M - 2\, \intr | M_{ij} |^2 \frac 1u \, dx &= -  2 \intr D_{ij} f \, D_i f \, D_j f \frac 1u \, dx \\&  -   \intr   |D f|^4 \, \frac 1u  dx -   \intr R \,  |D f|^2 \, dx
+ J^2_\rho.
\end{split}
\end{equation}
where
$$J^2_\rho = \int_{\partial B_\rho} R  \, \frac{\partial f}{\partial n} \, d\sigma -   \int_{\partial B_\rho}  \frac{\partial (|Df|^2)}{\partial n}  \, \frac 1u \, d\sigma.$$
We next observe that 
$$2 \intr D_{ij} f \, D_i f \, D_j f \frac 1u \, dx = \intr D_i (|D f|^2) \, D_i f \frac 1u$$
and integrate  by parts using once more that $\Delta f=R\, u$ and that 
$D_i f = - u^{-1} D_i f$, to find
$$
2 \intr D_{ij} f \, D_i f \, D_j f \frac 1u \, dx =  -   \intr R\,  |D f|^2 \, dx - 
\intr   |D f|^4 \, \frac 1u  dx   + J_\rho^3
$$
where 
$$J_\rho^3=  \lim_{\rho \to \infty} \int_{\partial B_\rho} |Df|^2 \, \frac{\partial f}{\partial n} \, d\sigma.$$
Combining the above we conclude that 
$$ M - 2\, \intr | M_{ij} |^2 \frac 1u \, dx =  J_\rho$$
with
$$J_\rho = \int_{\partial B_\rho} R  \, \frac{\partial f}{\partial n} \, d\sigma  -  \int_{\partial B_\rho}  \left ( \frac{\partial (|Df|^2)}{\partial n}  +   |Df|^2 \, \frac{\partial f}{\partial n} \right ) \, \frac 1u  \, d\sigma. 
$$
We will now show that $\lim_{\rho \to \infty} J_\rho =0.$ Clearly 
the first term tends to zero, because 
$\int_{\partial B_\rho} |{\partial f}/{\partial n}| \, d\sigma$ is  bounded Lemma \ref{lem-3}  and  $R(x,t) \to 0$, as $|x| \to \infty$, by Lemma \ref{lem-4}.   
It remains to show that
\begin{equation}\label{eqn-vbt}
\lim_{\rho \to \infty} \int_{\partial B_\rho}  \left ( \frac{\partial (|Df|^2)}{\partial n}  +   |Df|^2 \, \frac{\partial f}{\partial n} \right ) \, \frac 1u  \, d\sigma =0.
\end{equation}

We first observe that since $f = -\log u$, we have
$$\left ( \frac{\partial (|Df|^2)}{\partial n}  +   |Df|^2 \, \frac{\partial f}{\partial n} \right ) \frac 1u
=  \frac{\partial}{\partial n} \left ( \frac{|D\log u|^2}{u} \right ).$$
Expressing the last  term in  cylindrical coordinates, setting  
$v(s,\theta,t) = r^2 \, u(r,\theta,t)$,  with $s=\log r$,  we find
$$ \int_{\partial B_\rho}  \frac{\partial}{\partial n} \left ( \frac{|D\log u|^2}{u} \right ) \frac 1u \, d\sigma  = \int_0^{2\pi}  \frac{\partial}{\partial s} \left ( 
\frac{[2-(\log v)_s]^2 + [(\log v)_\theta]^2}{v}  \right )\, d\theta.$$
Further computation shows that
\begin{equation*}
\begin{split}
\frac {\partial}{\partial s} & \left ( 
\frac{[2-(\log v)_s]^2 + [(\log v)_\theta]^2}{v}  \right )\,  =
- 2\, [2 - (\log v)_s] \, \frac{(\log v)_{ss}}v \\&+ 2\, (\log v)_\theta \, \frac{(\log v)_{\theta s}}{v} + \{ [2-(\log v)_s]^2 + [(\log v)_\theta]^2 \} \, (\log v)_s. 
\end{split}
\end{equation*}
By Lemma \ref{lem-3} , $(\log v)_\theta$ is bounded as $s \to \infty$, while
$(\log v)_{s}$, $(\log v)_{\theta s}$ and $(\log v)_{ss}$ tend to
zero, as $s \to \infty$. Using also that $v$ is bounded away from
zero as $s \to \infty$, we finally conclude that
$$\lim_{s \to \infty} \int_0^{2\pi}  \frac{\partial}{\partial s} \left ( 
\frac{[2-(\log v)_s]^2 + [(\log v)_\theta]^2}{v}  \right )\, d\theta =0$$
implying \eqref{eqn-vbt} therefore finishing the proof of the claim
\eqref{eqn-eter33}.

We will now  conclude the proof of the Theorem. From \eqref{eqn-eter3} and \eqref{eqn-eter33} it follows that
$$\intr \frac{|X|^2}R \, dx + 2\, \intr | M_{ij} |^2 \frac 1u \, dx \leq I_\rho 
+ J_\rho$$
where both $$\lim_{\rho \to \infty} I_\rho + J_\rho = 0.$$
This immediately gives that $X\equiv 0$ and $M_{ij} \equiv 0$
for all $t$ showing that $U$ is a gradient soliton. It has been
shown by L.F. Wu \cite{W2} that there are only two types of
complete gradient solitons on $\R^2$,  the standard flat metric ($R\equiv 0$)
which is stationary,  and the cigar solitons \eqref{eqn-u}.    The flat solitons
violate condition \eqref{eqn-pw}. Hence, $u$ must be  of the form \eqref{eqn-soliton}, finishing the proof of the Theorem.  
\qed

\end{document}